\documentclass[11pt]{article}
\usepackage{amssymb,amsfonts,amsmath,latexsym,epsf,tikz,url}

\newtheorem{theorem}{Theorem}[section]

\newtheorem{corollary}[theorem]{Corollary}

\newtheorem{definition}[theorem]{Definition}

\newtheorem{problem}[theorem]{Problem}

\newcommand{\proof}{\noindent{\bf Proof. }}
\newcommand{\qed}{\hfill $\square$\medskip}

\begin{document}

\title{Trees  with  distinguishing number two}

\author{
Saeid Alikhani  $^{}$\footnote{Corresponding author} \and Samaneh Soltani
}

\date{\today}

\maketitle

\begin{center}
Department of Mathematics, Yazd University, 89195-741, Yazd, Iran\\
{\tt  alikhani@yazd.ac.ir,  s.soltani1979@gmail.com}
\end{center}

%%%%%%%%%%%%%%ABSTRACT%%%%%%%%%%%%%%%%%%%%%%%%%%%%%%%%%%%%%%%%%%%%%%%%%%%%%%%%%%%%

\begin{abstract}
The distinguishing number $D(G)$  of a graph $G$ is the least integer $d$
such that $G$ has a vertex labeling   with $d$ labels  that is preserved only by a trivial
automorphism.  In this paper we characterize all trees with radius at most three  and distinguishing number two. Also we present a necessary condition for trees with distinguishing number two and radius more than three. 
\end{abstract}


\medskip
\noindent{\bf AMS Subj.\ Class.}: 05C25 

%%%%%%%%%%%%%%%%%%%%%%%%%%%%%%%%%%%%%%%%%%%%%%%%%%%%%%%%%%%%%%%%%%%%%%%%%%%%%%%%%
%%%%%%%%%%%%%%%%%%%%%%%%%%%%%%%%%%%%%%%%%%%%%%%%%%%%%%%%%%%%%%%%%%%%%%%%%%%%%%%%%
\section{Introduction and definitions}
%%%%%%%%%%%%%%%%%%%%%%%%%%%%%%%%%%%%%%%%%%%%%%%%%%%%%%%%%%%%%%%%%%%%%%%%%%%%%%%%%
Let $G=(V,E)$ be a simple connected graph. The automorphism group and the complement of $G$, is denoted by  ${\rm Aut}(G)$  and $\overline{G}$, respectively.  The \textit{distance}   between two vertices $u, v$ of $G$, denoted by $d_G(u, v)$, is
the length of a shortest $u - v$ path if any; otherwise, $d_G(u, v) = \infty$. The  \textit{eccentricity } $e_G(v)$ of a vertex $v$ in $G$ is the distance
between $v$ and a vertex farthest from $v$. A vertex $u$ is an \textit{eccentric vertex of a vertex}  
$v$ in $G$ if $d_G(v, u) = e_G(v)$. The minimum eccentricity among the vertices of $G$ is the 
 \textit{radius} of $G$, ${\rm rad}(G)$, while the maximum eccentricity is its \textit{diameter}  ${\rm diam}(G)$. A vertex
$v$ of $G$ is called a \textit{central vertex}  if $e_G(v) = {\rm rad}(G)$, and center $C(G)$ of $G$ is the
collection of all such vertices in $G$. Graph $G$ with ${\rm rad}(G) = {\rm diam}(G)$ is  called
 \textit{self-centered} graph. Equivalently, a graph is selfcentered if all its vertices lie in the center. Thus, the eccentric set of a
self-centered graph contains only one element, that is, all the vertices have the same
eccentricity. It follows that any disconnected graph is self-centered graph
with radius $\infty$. The subgraph induced by center of $G$, $C(G)$, in $G$ is denoted by
${<{C(G)}>}_G$. A graph $G$ is called  \textit{$k$-self-centered} if ${\rm diam}(G) = {\rm rad}(G) = k$. The terminology \textit{$k$-equi-eccentric} graph  is also used by some authors. For studies on these graphs see \cite{F. Buckley,A.C. Malaravan, S. Negami and G.H Xu} and \cite{M.H. Shekarriz}.

Distinguishing labeling was first defined by Albertson and Collins \cite{Albert} for graphs. A labeling of a graph $G$, $\phi : V(G)\rightarrow \{1,2,\ldots ,r\}$, is said to be \textit{$r$-distinguishing} if no nontrivial automorphism of $G$ preserves all the vertex labels. In other words, $\phi$ is $r$-distinguishing if for any $\sigma\in {\rm Aut}(G)$, $\sigma \neq id$, there is a vertex $x$ such that $\phi (x) \neq \phi (\sigma (x))$. The \textit{distinguishing number} of a graph $G$ is defined as 
\begin{equation*}
D(G) = {\rm min}\{r \vert ~ G ~\textsl{\rm{has a labeling that is $r$-distinguishing}}\}.
\end{equation*} 

It is immediate that $D(K_n) = n$ for the complete graph $K_n$ of order $n$, and that $D(P_n) = 2$ for $n\geq 2$, where $P_n$ is the path graph of order $n$.   Cheng in \cite{C.T. Cheng} has presented $O(nlog n)$-time algorithms that compute the distinguishing numbers of trees
and forests. The distinguishing number of graphs have been extensively studied in the
literature.  
 The following general problem  appeared in \cite{Open problems column}:
 
 \begin{problem}
 Characterize graphs with distinguishing number two. 
  \end{problem} 
 
 In this paper,  we obtain all trees with radius at most two and distinguishing number two.  Also we state a note on  the trees of radius  greater than two and finally pose a problem.

\section{ Main results}

To obtain all trees with radius at most two and distinguishing number two, we need some  preliminaries:

\begin{definition}
 Let $G$ be a 2-self-centered graph. A vertex $x$ in $G$ is called critical
for $u$ and $v$ if $uv \notin E$ and $x$ is the only common neighbour of $u$ and $v$.
\end{definition}
\begin{definition}
Let $G$ be a graph. A cycle $C$ in $G$ is said to be locally geodesic at a vertex $v$ if for each vertex $u$ on
$C$, the distance between $v$ and $u$ in $C$ coincides with that in $G$.
\end{definition}
\begin{corollary} {\rm \cite{S. Negami and G.H Xu}}\label{Corollary7} Let $G$ be a bipartite graph. Then the following three are equivalent:
\begin{itemize}
\item[(i)] $G$ is isomorphic to the complete bipartite graph $K_{n,m}$ $(n, m \geqslant 2)$;
\item[(ii)] $G$ is $2$-self-centered;
\item[(iii)] $G$ is a block and for each vertex $v$ of $G$, there is no cycle locally geodesic at
$v$ and of length more than $4$.
\end{itemize}
\end{corollary}

Now, we can consider edge-minimal 2-self-centered graphs with some triangles.  We need the following procedure to proceed.

\medskip

\textbf{Procedure.}{\rm \cite{M.H. Shekarriz}} 
Let $G$ be a graph, $u, v,w$ form a triangle in $G$ and suppose that $v$ is a critical vertex for $u$ and $v_1,\ldots , v_p$ and/or $u$ is a critical vertex for $v$ and
$u_1,\ldots , u_q$. Remove the edge $uv$ and add edges $uv_1,\ldots , uv_p$ and $vu_1,\ldots , vu_q$. The following theorem characterizes edge-minimal $2$-self-centered graphs with
 triangles, which completes the characterization of all 2-self-centered graphs.
\begin{theorem}{\rm \cite{M.H. Shekarriz}}\label{Theorem15} Let $G$ be a graph. Then $G$ is an edge-minimal $2$-self-centered graph with some triangle if and only if the following two conditions are true:
\begin{itemize}
\item[(i)] For each edge of every triangle in $G$, at least one end-vertex is a critical
 vertex (for the other end-vertex of that edge and some other vertices of $G$), and
\item[(ii)] Iteration of above Procedure  on $G$ (at most to the number of triangles of $G$)
 transforms $G$ to a triangle-free $2$-self-centered graph.
\end{itemize}
\end{theorem}

The sequential join $G_1 + G_2 + \cdots + G_t$, of graphs $G_1, G_2,\ldots , G_t$, is formed from
$G_1\cup G_2\cup \ldots \cup G_t$, by adding the additional edges joining each vertex of $G_k$ with
each vertex of $G_{k+1}$, $1 \leqslant k \leqslant t - 1$. The double star $S_{a,b}$ is the tree $\overline{K}_a + K_1+K_1 + \overline{K}_b$, and the tritip $K_3(a, b, c)$ is the
graph formed by adding $a, b$, and $c$ pendant edges at the vertices of $K_3$, see  \cite{F. Buckley}.

\begin{figure}
	\begin{center}
\includegraphics[width=0.6\textwidth]{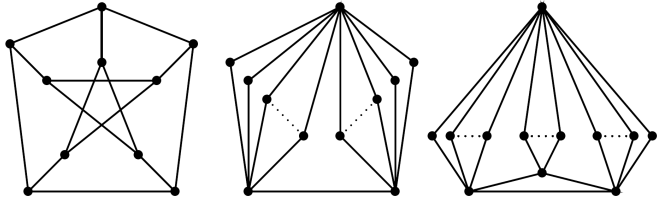}
	\caption{\label{fig1} The edge-minimal $2$-self-centered graphs in Theorem \ref{Theorem7}.}
\end{center}	
\end{figure}

\begin{theorem}{\rm \cite{F. Buckley}}\label{Theorem7}. Let $G$ be a self-centered graph with $n \geqslant 5$ vertices and diameter two,
having as few edges as possible among such graphs. Then $G$ is one of the following (see
Figure \ref{fig1}):
\begin{itemize}
\item[(i)] The Petersen graph,
\item[(ii)] The graph formed from $S_{a,b}$ by adding an additional vertex $u$ and joining $u$ to each vertex of degree 1 in $S_{a,b}$.
\item[(iii)] The graph formed from $K_3(a, b, c)$ by adding a new vertex $w$ and joining $w$ to
each vertex of degree 1 in $K_3(a, b, c)$, $a + b + c = n - 4$, $a, b, c \geqslant 1$.
\end{itemize}
\end{theorem}

\begin{theorem}{\rm  \cite{A.C. Malaravan}}\label{Theorem3.3}
	Let $T$ be a tree. Then, ${<C(\overline{T})>}_{\overline{T}}\cong \overline{T}$ if and only if ${\rm diam}(T ) \neq 3$.
\end{theorem}
\begin{corollary}{\rm  \cite{A.C. Malaravan}}\label{Corollary3.4}
Let $T$ be a tree of order $n$. Then, ${<C(\overline{T})>}_{\overline{T}}$ is either $\overline{T}$ or isomorphic to $K_{n-2}$.
\end{corollary}

\begin{figure}
	\begin{center}
		\includegraphics[width=1.0\textwidth]{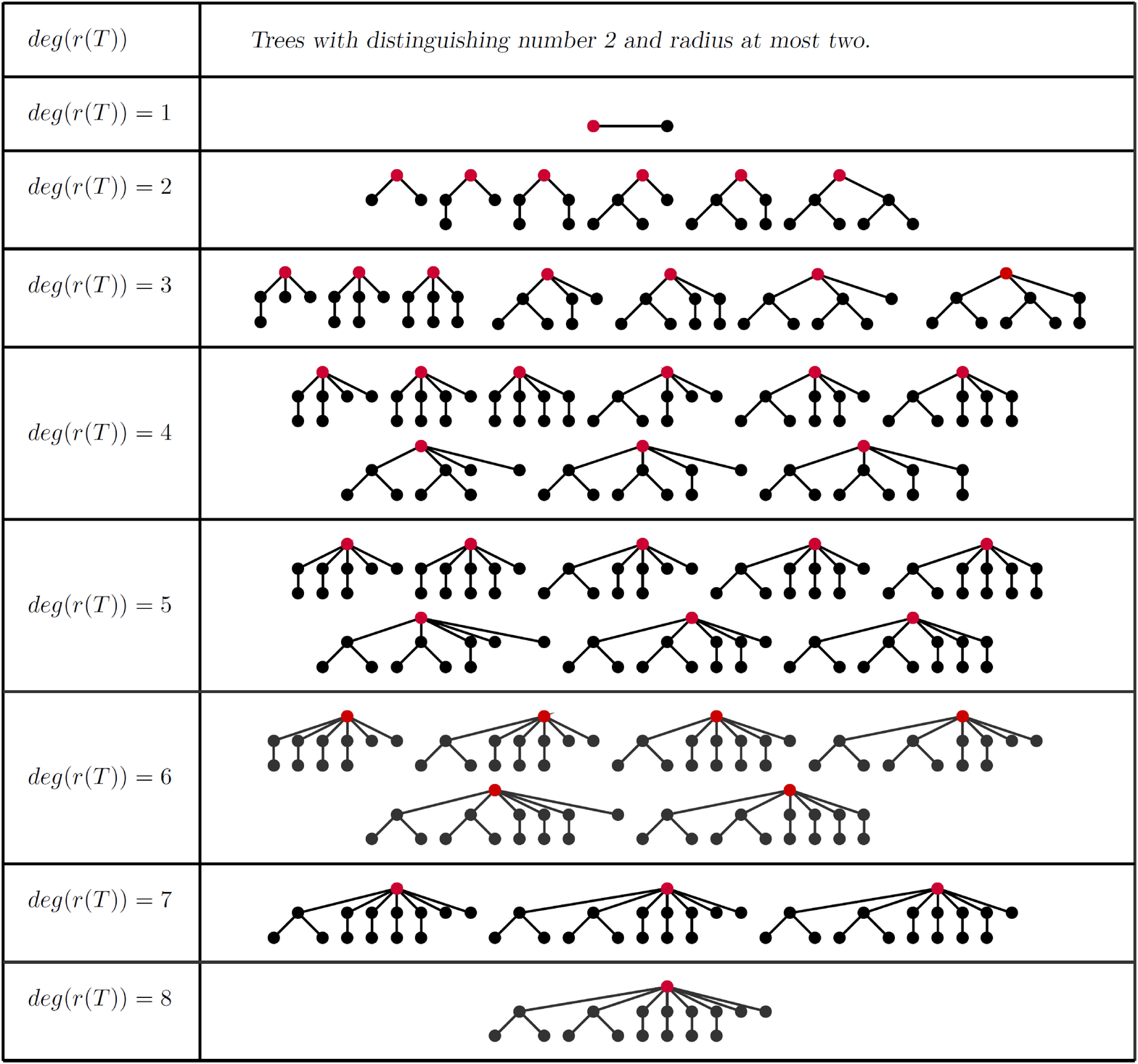}
		\caption{\label{fig0} Trees  with distinguishing number $2$ and radius at most two.}
	\end{center}	
\end{figure}

\begin{theorem}
Let  $T$ be a tree of radius at most two. The  distinguishing number of $T$ is two if and only if $T$ is one of the trees in Figure \ref{fig0}.
\end{theorem}
\proof We know  that given an unrooted tree $T'$ we can construct a rooted tree $T$ such that
$D(T') = D(T)$ \cite{C.T. Cheng}.   It is well known that a tree either has one center
(i.e., it is unicentral) or has two adjacent centers (i.e., it is bicentral). Thus, if $T'$ has a unique center, simply let $T$ be a copy of $T'$;
otherwise, let $T$ be the tree formed by appending a new vertex to the two centers of $T'$ and
deleting the edge between the two old centers of $T'$. In both cases, $T$ has a unique center which we designate as its root $r(T)$. So for obtaining all trees with radius at most two and the distinguishing number two it is sufficient to consider all rooted trees of radius at most two with distinguishing number two.  Only with respect to the degree of the root, we can get these trees as shown in Figure \ref{fig0}. Note that the roots have shown in red colour. \qed

\medskip
Here, we state a note about trees with the distinguishing number two and radius more than two.  First we state and prove the following theorem: 

\begin{theorem} \label{thmrad}
There is no tree $T$ with diameter tree, radius ${\rm rad}(T)\geqslant 3$ and $D(T)=2$.  	
	\end{theorem} 
\proof 
Since  ${\rm diam}(T)=3$, so  the subgraph induced by center of $\overline{T}$, $C(\overline{T})$, in $\overline{T}$, is ${<C(\overline{T})>}_{\overline{T}}\cong K_{n-2}$ by Corollary \ref{Corollary3.4}.  So the only possible number for $n$ for which $D(T)=2$ is $n=4$, and hence $T=P_4$. But ${\rm rad}(T)\geqslant 3$, and therefore there is no tree with distinguishing number two in this case.\qed

\begin{theorem} \label{geq4}
If $T$ is a tree with radius ${\rm rad}(T) \geqslant 3$ and $D(T)=2$, then	$\overline{T}$  is a connected $2$-self-centered graph that is not bipartite.
\end{theorem} 
\proof	
If $T$ is a tree of order $n$ with radius ${\rm rad}(T) \geqslant 3$ and $D(T)=2$, then the radius of the complement of $T$ is ${\rm rad}(\overline{T})\leqslant 2$. By Theorem \ref{thmrad}, ${\rm diam}(T)\geqslant 4$, and then ${<C(\overline{T})>}_{\overline{T}}\cong \overline{T}$ by Theorem \ref{Theorem3.3}. It is known that the complement of a tree is connected or it is a union of an isolated vertex and a complete graph. If $\overline{T}$ is a union of an isolated vertex and a complete graph, then $T$ is a star graph, and since $D(T)=2$ so $T=P_3$. But ${\rm diam}(T)\geqslant 4$, thus there is no tree with distinguishing number two in this case. Therefore $\overline{T}$  is a connected $2$-self-centered graph that is not bipartite by Corollary \ref{Corollary7} and Theorem \ref{Theorem3.3}. \qed

\medskip 

By Theorem \ref{geq4},  for finding all trees with ${\rm rad} (T)\geqslant 3$ and $D(T)=2$, we  should only consider the $2$-self-centered trees that are not bipartite. But Buckley in \cite{F. Buckley} obtained all $2$-self-centered graphs having as few edges as possible among such graphs. Since $T$ is a tree, so $\overline{T}$ is not an edge-minimal $2$-self-centered graph by Theorem \ref{Theorem15}. Thus we should add a specific  number of edges to each edge-minimal $2$-self-centered graph so that we get $\overline{T}$. We end this paper with the following problem: 

\begin{problem} 
	Characterize edges which  add  to each  edge-minimal $2$-self-centered graph such  that $D(\overline{T})=2$.
\end{problem}

\end{document}